\documentclass[12pt]{article}
\usepackage[english]{babel}
\usepackage{latexsym,amsfonts,amssymb,amsmath,enumerate,amsthm}

\title{Second derivatives of norms and contractive
complementation in vector-valued spaces}

\author{Bas Lemmens, Beata Randrianantoanina\footnote{Participant, NSF Workshop
in Linear Analysis and Probability, Texas A\&M University.}, \\
and Onno van Gaans}

\date{}
\let\epsilon=\varepsilon
\let\phi=\varphi
\let\theta=\vartheta
\newcommand{\NN}{\mathbb{N}}
\newcommand{\RR}{\mathbb{R}}
\newcommand{\linspan}{\mathrm{span}\,}

\newtheorem{theorem}{Theorem}[section]
\newtheorem{definition}[theorem]{Definition}
\newtheorem{lemma}[theorem]{Lemma}
\newtheorem{corollary}[theorem]{Corollary}
\newtheorem{proposition}[theorem]{Proposition}

\renewcommand{\descriptionlabel}[1]%
    {\hspace{\labelsep}\textrm{#1}}

\begin{document}
\maketitle
\begin{abstract}
We consider 1-complemented subspaces (ranges of contractive
projections) of vector-valued spaces $\ell_p(X)$, where $X$ is a
Banach space with a 1-unconditional basis and $p\in (1,2)\cup
(2,\infty)$. If the norm of $X$ is twice continuously differentiable
and satisfies certain conditions connecting the norm and the notion
of disjointness with respect to the basis, then we prove that every
1-complemented subspace of $\ell_p(X)$ admits a basis of mutually
disjoint elements. Moreover, we show that every contractive
projection is then an averaging operator. We apply our results to
the space $\ell_p(\ell_q)$ with $p,q\in (1,2)\cup (2,\infty)$ and obtain 
a complete characterization of its 1-complemented
subspaces.
\end{abstract}
\vspace{4mm}
\begin{small}
\noindent
\textbf{M.S.C. (2000)}. Primary 46B45, 46B04; Secondary 47B37.
\end{small}
\newpage
\section{Introduction}\label{sec:1}
It is well-known that every orthogonal projection on a Hilbert space $H$ is
contractive and that for every closed subspace $Y$ of $H$ there exists an orthogonal
projection whose range is $Y$.
In fact, if the dimension of $H$ is at least 3, then it has been shown by Kakutani
\cite{Ka} that this property characterizes the Hilbert spaces among the Banach spaces.
Hence contractive projections on Banach spaces are a natural generalization of
orthogonal projections on Hilbert spaces; but, as the result of Kakutani indicates,
they are more rare.
It is therefore not surprising that contractive projections have been studied
extensively, starting with work by Bohnenblust \cite{Bo} in the nineteen
forties.
A detailed overview of the vast literature on contractive projections and their applications 
can be found in the survey papers \cite{CP,Ran5}.

Identifying contractive projections and their ranges for a given Banach space has
often proved to be difficult.
For Lebesgue $L_p$ spaces, with $p\neq 2$, there exists a well-known characterization of the 
contractive projections and their ranges (see \cite{An,BL,Do,LT,Tz}). 
However, for many other classical Banach spaces, such as Orlicz spaces and Lorentz sequence 
spaces, there are only partial results; see the recent survey \cite{Ran5}. 
The best results to date for Musielak-Orlicz spaces were obtained recently by Jamison, Kami\'nska, 
and Lewicki \cite{JKL}, who characterized the 1-complemented subspaces of finite co-dimension 
in case the Orlicz functions are sufficiently smooth.   

The problem of characterizing contractive projections and their
ranges in vector-valued spaces is known to be particularly hard,
even for spaces of the type $L_p(X)$, where $X$ is finite
dimensional. We refer the reader to a survey by Doust \cite{Doust} for
a nice overview of the various partial results that were known in the 
mid nineties. As Doust points out, most of the results require special 
additional assumptions on the form of projections
and leave open the cases of $L_p(X)$ and $\ell_p(X)$, even when $X$ is a two 
dimensional $\ell_q$ space.
The most general result for vector-valued spaces was obtained recently by
Raynaud \cite{Ray}, who gave a complete description of the contractive
projections in $L_p(H)$, for $H$ a Hilbert space. Raynaud's result is
valid in full generality without any assumptions about $\sigma-$finiteness
of the measure on $L_p$ or separability of Hilbert space $H$.

In the present paper we study contractive projections in
vector-valued $\ell_p(X)$ spaces, in particular $\ell_p(\ell_q)$ spaces. 
Our main result characterizes the 1-complemented subspaces of real 
$\ell_p(\ell_q)$ spaces for $p,q \in (1,2)\cup(2,\infty)$.
The complex case was obtained in \cite{Ran1}, in which 
1-complemented subspaces of general complex spaces with 1-unconditional bases 
are characterized, using a technique based on Hermitian operators. 
This technique, however, does not work in real spaces. 
This is due to the following fundamental result of Kalton and Wood \cite{KW} for complex spaces, 
which has no real analogue (see \cite{Le79}): 
every 1-complemented subspace of a complex Banach space with a 1-unconditional basis 
has a 1-unconditional basis.  
In fact, it is not known whether each 1-complemented subspace in real space admits an 
unconditional basis with any constant (cf. \cite[Section 7.e]{Ran5}).
The characterization presented here shows in particular that every 1-complemented subspace 
of areal $\ell_p(\ell_q)$ space has an unconditional basis.

To obtain the characterization, we introduce a condition on the second derivative 
of the norm of a real Banach space, with a $1$-unconditional basis, that guarantees that 
every contractive projection is an averaging operator and its range admits a block
basis. 
We subsequently apply it to analyse contractive projections on
vector-valued $\ell_p(X)$ spaces. 
As a consequence we find that 1-complemented subspaces of 
$\ell_p(\ell_q)$ admit a block basis, in case $p$ and $q$ both in $(2,\infty)$, 
or, both in $(1,2)$.
In the mixed case, where $p\in (1,2)$ and $q\in (2,\infty)$, or, the other way around,
we can not apply the condition and we shall use a different argument.

The idea to exploit the second derivative of the norm to analyse contractive
projections is not new. In fact, it already appears implicitly in work of
Bernau and Lacey \cite{BL} and Lindenstrauss and Tzafriri \cite{LT}, who
considered the derivative of the norming functionals (or duality map).
The derivative of the norming functionals has also been used by Bru, Heinich, and
Lootgieter \cite{BHL} to identify contractive projections on Orlicz spaces that
have a second order smooth norm and satisfy some additional constraints.
More recently, Lemmens and van Gaans \cite{LvG} have used the second derivative of
the norm to show for a fairly general class of finite dimensional Banach spaces
that the range of every contractive projection has a block basis.
In particular, one could apply their technique to prove that the range of each
contractive projection on $\ell_p^m(\ell_q^n)$, where $m,n\in\mathbb{N}$, has
a block basis, if $p>2$ and $q>2$, or $p$ and $q$ in $(1,2)$. 
It is interesting to note that differential techniques are not only useful for analysing 
contractive projections, but also appear in the study of isometries; see, for instance, 
work by Koldobsky \cite{Kol}.

\section{Preliminaries}\label{sec:2}
In this section we collect several definitions and notations that will be
used  throughout the exposition.
In addition, we recall some preliminary results.
Before we get started however, we emphasize that in the paper all Banach spaces are
over the field of real numbers.

Let $X$ be a Banach space with a $1$-unconditional basis $\{e_i\}_{i=1}^\infty$.
We denote by $S_X=\{x\in X\colon \|x\|=1\}$ the \emph{unit sphere of} $X$.
The dual space of $X$ is denoted by $X^*$ and the conjugate norm is indicated by
$\|\cdot\|^*$.
For each $x=\sum_{i=1}^\infty \alpha_ie_i$ in $X$ we let
$s(x)=\{i\in \mathbb{N}\colon \alpha_i\neq 0\}$ be the \emph{support} of $x$.
If $s(x)$ is finite, we call $x$ \emph{simple}. 
Since $X$ has a $1$-unconditional basis, the partial ordering induced by the basis
makes $X$ a Banach lattice.
Hence there exists a natural notion of disjointness in $X$.
Elements $x$ and $y$ in $X$ are called \emph{disjoint} if $|x|\wedge |y|=0$.
As $X$ has a $1$-unconditional basis, this is equivalent to $s(x)\cap s(y)=\emptyset$.
We note that the dual of a Banach lattice is again a Banach lattice and therefore
there exists a natural notion of disjointness in $X^*$.
If $Y$ is a subspace of $X$ and $Y$ has a basis $\{y_i\}_{i=1}^{\dim Y}$ such that
$y_i$ and $y_j$ are disjoint for all $i\neq j$, then $\{y_i\}_{i=1}^{\dim Y}$ is
called a \emph{block basis} for $Y$.

In this paper we are particularly interested in the vector-valued spaces $\ell_p(X)$.
If $X$ is a Banach space with norm $\sigma$, then for $p\in [1,\infty)$ the vector-valued space
$\ell_p(X)$ consists of those $x\colon \mathbb{N}\to X$ for which
\[
\Big{(}\sum_{k=1}^\infty \sigma(x(k))^p\Big{)}^{1/p}< \infty.
\]
By equipping the space $\ell_p(X)$ with the norm,
\[
\|x\|_{\ell_p(X)}= \Big{(}\sum_{k=1}^\infty \sigma(x(k))^p\Big{)}^{1/p}
\mbox{\quad  for all } x\in \ell_p(X),
\]
it becomes a Banach space.
It is not difficult to verify that if $X$ has a $1$-unconditional basis
$\{e_i\}_{i=1}^\infty$, then for any ordering on the elements
$(i,j)\in \mathbb{N}\times \mathbb{N}$ the functions
$e_{ij}\colon \mathbb{N}\to X$ given by, $e_{ij}(k)=e_j$ if $k=j$ and
$e_{ij}(k)=0$ otherwise, form a $1$-unconditional basis for $\ell_p(X)$.
For elements in $\ell_p(X)$ with basis $\{e_{ij}\}_{(i,j)}$, it is useful to
introduce the notion of vector support.
For $x\in \ell_p(X)$ we define $vs(x)=\{k\in\mathbb{N}\colon x(k)\neq 0\}$ to be the
\emph{vector support} of $x$.
We note that the dual space of $\ell_p(X)$ is equal to $\ell_{p^*}(X^*)$, where
$1/p+1/p^*=1$, if $X$ is reflexive (see e.g. \cite[Chapter~IV]{DU}).

If $X$ is a Banach space with a $1$-unconditional basis and $T\colon X\to X$ is a
linear operator for which there exist mutually disjoint elements $\{u_i\}_{i\in I}$
in $X$ and mutually disjoint elements $\{v_i^*\}_{i\in I}$ in $X^*$ such that
$v_i^*(u_j)=0$ for all $i\neq j$ and
\[
Tx=\sum_{i\in I} v_i^*(x)u_i\mbox{\quad for all }x\in X,
\]
then $T$ is called an \emph{averaging operator}. 
Obviously, the range of an averaging operator has a block basis.

We also need to recall some definitions concerning higher order derivatives of norms.
Let $X$ and $Y$ be Banach spaces and let $L(X,Y)$ be the Banach space of continuous
linear operators from $X$ into $Y$ equipped with the usual operator norm.
We denote by $B^k(X,Y)$ the Banach space of continuous $k$-linear operators
$T\colon X\times \cdots\times X\to Y$ with the norm
$\|T\|=\sup\{ \|T(x_1,\ldots,x_k)\| \colon \|x_1\|=\ldots=\|x_k\|=1\}$.
A mapping $\phi\colon U\to Y$, where $U\subset X$ is open, is called
\emph{differentiable at} $x\in U$ if there exists a linear operator
$D\phi(x)$ in $L(X,Y)$ such that
\[
\lim_{h\to 0} \frac{\|\phi(x+h)-\phi(x) -D\phi(x)h \|}{\|h\|}=0.
\]
The linear operator $D\phi(x)$ is unique and is called the \emph{derivative} of
$\phi$ at $x$. Higher order derivatives $D^k\phi(x)\in B^k(X,Y)$ are defined in the
usual inductive manner (see Dieudonn\'e \cite{Dieu}).
The mapping $\phi\colon U\to Y$ is said to be $C^k$ on $U$ if it is $k$-times
differentiable at every point $x\in U$ and $D^k\phi\colon U\to B^k(X,Y)$ is continuous.
We note that if $\phi\colon U\to Y$ is $k$-times differentiable at $x$, then the
multi-linear map $D^k\phi(x)$ is symmetric.

If $\phi\colon U\to Y$ is $C^k$ on $U$ and $[x,x+h]\subset U$, then the usual Taylor
expansion
\[
\phi(x+h)=\phi(x)+\sum_{j=1}^k \frac{1}{j!}D^j\phi(x)(h,\ldots,h) + \theta_x(h),
\]
where
\[
\lim_{\|h\|\to 0} \frac{\|\theta_x(h)\|}{\|h\|^k}=0,
\]
is valid.
We shall also use a weaker notion of differentiability.
A map $\phi\colon U\to Y$ is said to be \emph{$k$-times directionally
differentiable at} $x\in U$ if for $1\leq j\leq k$ there exists a continuous
symmetric $j$-linear
operator $D^j\phi(x)\in B^j(X,Y)$ such that for every $h\in X$ and $t\in \mathbb{R}$
with $[x,x+th]\subset U$ we have that
\[
\phi(x+th)=\phi(x)+\sum_{j=1}^k \frac{t^j}{j!}D^j\phi(x)(h,\ldots,h) + \theta_x(th),
\]
where
\[
\lim_{t\to 0} \frac{\|\theta_x(th)\|}{|t|^k}=0.
\]
One can verify that the operators $D^j\phi(x)$ are unique.
Equipped with these notions of differentiability we now recall the following
definition from \cite{SS}.
\begin{definition}\label{def:2.1}
A Banach space $X$ is called \emph{$C^k$ smooth} if the norm $\|\cdot\|$ is $C^k$ on
$X\setminus\{0\}$.
It is said to be \emph{$D^k$ smooth} if the norm is $k$ times directionally
differentiable at each each $x\in X\setminus\{0\}$.
\end{definition}
In the analysis of the vector-valued spaces $\ell_p(X)$ we shall use the following
smoothness result of Leonard and Sundaresan.
\begin{theorem}[\cite{LS}, Theorem 3.3]\label{thm:2.2}
If $X$ is a Banach space and $p>k$, then $\ell_p(X)$ is $C^k$ smooth
if and only if $X$ is $C^k$ smooth and the $k$th derivative of the norm of $X$ is
uniformly bounded on the unit sphere in $X$.
\end{theorem}
In particular, it follows from this theorem that the space $\ell_p(\ell_q)$ is
$C^2$ smooth if $p$ and $q$ in $(2,\infty)$.

\section{Second derivatives of norms and contractive projections}\label{sec:3}

We begin this section by introducing the property of the second derivative of the norm.
Subsequently we explain how it can be used to analyse contractive projections and their
ranges.
\begin{definition}\label{def:3.1}
Suppose $X$ is a Banach space with basis and for each $x,y\in X$ the function
$N=N_{xy}\colon \mathbb{R}\to \mathbb{R}$ is given by $N(\alpha)=\|x+\alpha y\|$
for all $\alpha\in \mathbb{R}$.
We say that $X$ \emph{reflects disjointness} if for every
$x,y\in X$ with $x\not\in\linspan\{y\}$ the following conditions hold:
\begin{enumerate}[(i)]
\item the function $N$ is continuously differentiable (and then  $N''(\alpha)$ exists almost
everywhere by convexity);
\item if $x$ and $y$ are not disjoint and $N'(0)=0$, then
$N''(\alpha)$ does not converge to $0$, as $\alpha \to 0$ along any subset of
$[0,1]$ of full measure;
\item if $x$ and $y$ are disjoint and $y$ is simple, then $N'(0)=0$ and $N''(\alpha)$
converges to $0$, as $\alpha \to 0$ along a subset of $[0,1]$ of full measure.
\end{enumerate}
\end{definition}
The idea of this definition is that if $X$ reflects disjointness, then one can test
disjointness of elements in $X$ by analysing the second derivative of the norm.
Similar ways to test disjointness have been applied by Koldobsky \cite{Kol} to 
identify isometries.
The connection with contractive projections was found by Randrianantoanina
\cite{Ran6}, who used the condition to identify contractive projections
on certain Orlicz sequence spaces equipped with the Luxemburg norm.

Let us now  explain the connection of Definition \ref{def:3.1} with contractive
projections.
To do this it is useful to recall the following definition from \cite{Ran7}.
\begin{definition}\label{def:3.2}
Let $X$ be a Banach space with a basis and let $T\colon X\to X$ be a linear
operator. We call $T$ \emph{semi-band preserving} if for each $x,y\in X$ we
have that $Tx$ and $Ty$ are disjoint, whenever $Tx$ and $y$ are disjoint.
\end{definition}
Semi-band preserving operators on a Banach space with a $1$-unconditional basis
have a special form as the following theorem indicates.
\begin{theorem}[\cite{Ran7}, Theorem 4.7]\label{thm:3.3}
If $X$ is a Banach space and $X$ has a $1$-unconditional basis,
then a linear operator $T\colon X\to X$ is semi-band preserving if and only if
$T$ is an averaging operator.
\end{theorem}
Thus, to show that a contractive projection is an averaging operator and its
range has a block basis, it suffices to prove that the projection is semi-band
preserving.
Proving semi-band preservingness involves testing disjointness and this is where
the property in Definition \ref{def:3.1} comes into play.
As a matter of fact, we have the following result, which generalizes
\cite[Theorem 3.2]{Ran6}.
\begin{theorem}\label{thm:3.4}
If $X$ is a $D^2$ smooth Banach space, with a $1$-unconditional basis,
and $X$ reflects disjointness, then every contractive projection on $X$ is an
averaging operator and its range admits a block basis.
\end{theorem}
The proof of Theorem \ref{thm:3.4} is very similar to the proof of
\cite[Theorem 3.2]{Ran6} and uses the following lemma, which is a slight modification of \cite[Lemma 3.1]{Ran6}.
\begin{lemma}[\cite{Ran6}, Lemma 3.1]\label{lem:3.5}
Let $\phi\colon \mathbb{R}\to [0,\infty)$ and $\psi\colon \mathbb{R}\to [0,\infty)$
be convex functions that are continuously differentiable and assume that $\phi'$ and $\psi'$ are absolutely
continuous on $[0,1]$.
If $\phi(0)=\psi(0)$ and $\phi(\alpha)\leq \psi(\alpha)$ for all
$\alpha \in [0,\infty)$, then
\begin{enumerate}[(a)]
\item $\phi'(0)=\psi'(0)$;
\item the set $E=\{\alpha: \phi''(\alpha)$ and $\psi''(\alpha)$ exist and $\phi''(\alpha)\le \psi''(\alpha)\}$
has positive Lebesgue measure in each interval $(0,\delta)$, $\delta>0$;
\item for every $C>0$ the Lebesgue measure of
\[
\{\alpha\in [0,1]\colon \psi''(\alpha) \mbox{ exists and }\psi''(\alpha)\leq C\}
\]
is strictly smaller than $1$, whenever $\phi''(\alpha)\to \infty$\ as $\alpha\to 0$ along a subset of
full measure.
\end{enumerate}
\end{lemma}
\begin{proof}
The parts (a) and (c) are as in \cite[Lemma 3.1]{Ran6}. Part (b) is a modification of that lemma.
Let $A:=\{\alpha>0: \phi'(\alpha)=\psi'(\alpha)\}$. If $\inf A>0$, then there
is an $\varepsilon\in (0,\inf A)$ and $\phi'(\alpha)\neq \psi'(\alpha)$ for
all $\alpha\in (0,\varepsilon)$. Denote $h:=\psi-\phi$. Then $h\ge 0$ and $h'(\alpha)\neq 0$
for all $\alpha\in (0,\varepsilon)$. As $h'$ is continuous, $h'$ is either strictly positive
or strictly negative on $(0,\varepsilon)$. As $h(0)=0$ and $h\ge 0$, we have $h'(\alpha)>0$
for all $\alpha\in (0,\varepsilon)$. Since $h'$ is absolutely continuous, $h''$ exists almost
everywhere and
$$
h'(\alpha)=\int_0^\alpha h''(\beta)\,d\beta\quad\mbox{for all }\alpha\in (0,\varepsilon),
$$
as $h'(0)=0$. It follows that the set $E=\{\alpha\in (0,\varepsilon): h''(\alpha)\ge 0\}$ has
positive measure in each interval $(0,\delta)$, $\delta\in (0,\varepsilon)$.

If $\inf A=0$, then $\phi'(\alpha_n)=\psi'(\alpha_n)$ for some $\alpha_n\downarrow 0$. Then
$$
0=h'(\alpha_n)=\int_0^{\alpha_n} h''(\beta)d\beta
$$
and the conclusion follows.
\end{proof}

\begin{proof}[Proof of Theorem \ref{thm:3.4}]
By Theorem \ref{thm:3.3} it suffices to show that every contractive projection
$P\colon X\to X$ is semi-band preserving.
Let $x,y \in X$, with $y$ simple, and suppose that $Px$ and $y$ are disjoint and $Px\neq 0$.
Define functions
$\phi\colon \mathbb{R}\to [0,\infty)$ and $\psi\colon \mathbb{R}\to [0,\infty)$
by
\[
\phi(\alpha)=\|Px+\alpha Py\|\mbox{\quad and\quad }\psi(\alpha)=\|Px+\alpha y\|
\mbox{\quad for all }\alpha\in\mathbb{R}.
\]
Obviously, $\phi$ and $\psi$ are convex and $\phi(0)=\psi(0)$.
As $P$ is a contractive projection, $\phi(\alpha)=\|P^2x+\alpha Py\|\leq \psi(\alpha)$
for all $\alpha\in\mathbb{R}$.
Moreover, $\phi$ and $\psi$ are both twice continuously differentiable, because
$X$ is $D^2$ smooth.
We can now use the fact that $X$ reflects disjointness and $y$ is simple, to deduce
from Definition \ref{def:3.1}(iii) that $\psi'(0)=0$ and $\psi''(\alpha)$ converges
to $0$, as $\alpha\to 0$ along a subset of $[0,1]$ of full measure.
Since $\phi''$ is continuous,
Lemma \ref{lem:3.5} gives  $\phi'(0)=0$ and $\phi''(0)=0$.
By using Definition \ref{def:3.1}(ii) we find that $Px$ and $Py$ are disjoint. An arbitrary element $y\in X$
can be approximated by elements $y_n$ in $X$ such that the support of each $y_n$ is finite
and contained in the support of $y$.
We conclude that $P$ is semi-band preserving.
\end{proof}

The condition on the second derivative of the norm in Definition \ref{def:3.1}
has a natural interpretation in terms of curvature properties of the unit sphere,
if the Banach space is finite dimensional.
More precisely, one can show that if $X$ is a $C^2$ smooth finite dimensional Banach
space, with norm $\rho$ and a $1$-unconditional basis, then $X$ reflects disjointness
is equivalent to saying that for each $x\in S_X$ the normal curvature at $x$ in the
direction of $y$ is $0$ if and only if $x$ and $y$ are disjoint.
To prove this one has to note that the normal curvature $k(y)$ at $x$ in the direction 
$y$, where $y$ is in the tangent space at $x$, is given by 
$k(y)=N''_{xy}(0)/\|\nabla\rho(x)\|_2$.

\section{Contractive projections on $\ell_p(X)$} \label{sec:4}
In this section we analyse contractive projections on the vector-valued spaces
$\ell_p(X)$ and, in particular, the spaces $\ell_p(\ell_q)$, where $p$ and $q$
are not equal to $2$.
As mentioned in the introduction we distinguish two cases: the unmixed case, where
$p,q \in (2,\infty)$ or $p,q\in (1,2)$, and the mixed case, where $p\in (1,2)$ and
$q\in (2,\infty)$, or, the other way around.
We first prove a lemma and subsequently discuss the unmixed case.

\begin{lemma}\label{lem.L1}
Let $X$ be a $C^2$-smooth Banach space with norm $\sigma$ such that $D\sigma$ uniformly bounded 
on $S_X$. Let $p\in (1,2)\cup (2,\infty)$,
and let $x,y\in \ell_p(X)$. For $\alpha\in\RR$ define
\begin{eqnarray*}
N_k(\alpha)&:=&\sigma(x(k)+\alpha y(k)),\quad k\in\NN,\\
\tau(\alpha)&:=& \sum_{k=1}^\infty N_k(\alpha)^p,\\
N(\alpha)&:=&\tau(\alpha)^{1/p}.
\end{eqnarray*}
Then $\tau$ is $C^1$, $\tau'$ is absolutely continuous, $N'$ and $N''$ exist almost everywhere, and
\begin{eqnarray}
\tau'(\alpha) &=&
\sum_{k\in vs(x)\cap vs(y)} p N_k^{p-1}(\alpha)N_k'(\alpha)+p\alpha^{p-1}\sum_{k\in vs(y)\setminus vs(x)}
\sigma^p(y(k)),\label{f.tauprime}\\
\tau''(\alpha)&=& \sum_{k\in vs(x)\cap vs(y)} \Big( p(p-1)N_k^{p-2}(\alpha) N_k'(\alpha)^2 +
pN_k^{p-1}(\alpha)N_k''(\alpha)\Big)\nonumber\\
&&\qquad +p(p-1) \alpha^{p-2} \sum_{k\in vs(y)\setminus vs(x)} \sigma^p(y(k)),\label{f.taudouble}\\
\tau'(\alpha)&=&pN^{p-1}(\alpha)N'(\alpha), \label{f.tauprimeN}\\
\tau''(\alpha)&=& p(p-1)N^{p-2}(\alpha) N'(\alpha)^2+pN^{p-1}(\alpha)N''(\alpha)\label{f.taudoubleN}
\end{eqnarray}
for Lebesgue-almost every $\alpha\in\RR$.
\end{lemma}
\begin{proof}
Consider first a $k\in\NN$ such that $x(k)+\alpha y(k)\neq 0$ for all $\alpha\in\RR$.
Then $N_k$ is a $C^2$-function and hence $(N_k^p)'=pN_k^{p-1} N_k'$ is $C^1$.
Twice indefinite integration of $(N_k^p)''$ yields that
\begin{align}
&N_k^p(\alpha)=N_k^p(\beta) +pN_k^{p-1}(\beta) N_k'(\beta)(\alpha-\beta) \label{f.nkp}\\
&\qquad+\int_\beta^\alpha\int_\beta^s \Big( p(p-1)N_k^{p-2}(r)N_k'(r)^2 + pN_k(r)^{p-1}(r)N_k''(r)\Big)\,drds\nonumber
\end{align}
for every $\alpha,\beta\in\RR$, $\alpha\ge\beta$. For $k\in\NN$ such that $x(k)+\gamma y(k)= 0$
for some $\gamma\in\RR$, $N_k$ is $C^2$ on $\RR\setminus\{\gamma\}$ and a straightforward
computation shows that (\ref{f.nkp}) is true also in this case.

Next we integrate the right hand side of (\ref{f.taudouble}) twice.
As $N_k$ is convex, we have $N_k''\ge 0$ almost everywhere, so the right hand side of (\ref{f.taudouble})
is a measurable function with values in $[0,\infty]$. Fubini's theorem for positive functions and 
(\ref{f.nkp}) therefore yield that
\begin{align}
&\int_\beta^\alpha \int_\beta^s\left( \sum_{k\in vs(x)\cap vs(y)}\Big( p(p-1)N_k^{p-2}(r)N_k'(r)^2
+ pN_k^{p-1}(r)N_k''(r)\Big)\right. \nonumber\\
&\qquad\qquad\left.+ p(p-1)r^{p-2}\sum_{k\in vs(y)\setminus vs(x)} \sigma^p(y(k))\right)\,drds\nonumber\\
&\qquad= \tau(\alpha) -\tau(\beta) -\left( \sum_{k\in vs(x)\cap vs(y)} pN_k^{p-1}(\beta)
N_k'(\beta)\right)(\alpha-\beta)\nonumber\\
&\qquad\qquad -p\beta^{p-1} \sum_{k\in vs(y)\setminus vs(x)} \sigma^p(y(k))(\alpha-\beta),\label{f.twiceintegrated}
\end{align}
for every $\alpha\ge \beta$.
Since $D\sigma$ is uniformly bounded on $S_X$ and $D\sigma(x)=D\sigma(\lambda x)$ for all $\lambda \neq 0$, 
there exists a constant $C$ such that $\|D\sigma(z)\|\le C$
for all $z\in X\setminus\{0\}$. Due to Young's inequality we have
\begin{align}
|pN_k^{p-1}(\alpha)N_k'(\alpha)| &\leq pC\sigma(x(k))+\alpha y (k))^{p-1}\sigma(y(k))\nonumber\\
 \mbox{}&\le (p-1)C\Big(\sigma(x(k))+\sigma(y(k))\Big)^p + {\textstyle\frac{1}{p}}C\sigma(y(k))^p,\label{f.youngforprime}
\end{align}
and it therefore follows from (\ref{f.twiceintegrated}) that the right hand side of (\ref{f.taudouble}) is
an integrable function of $\alpha$ on bounded intervals.
By solving $\tau(\alpha)$ from (\ref{f.twiceintegrated}), it is clear that $\tau$ is $C^1$ and that
\begin{align*}
\tau'(\alpha)&=\sum_{k\in vs(x)\cap vs(y)} pN_k^{p-1}(\beta) N_k'(\beta) +
p\beta^{p-1}\sum_{k\in vs(y)\setminus vs(x)} \sigma^p(y(k))\\
&\qquad +\int_\beta^\alpha \left( \sum_{k\in vs(x)\cap vs(y)}
\Big( p(p-1)N_k^{p-2}(r)N_k'(r)^2+pN_k^{p-1}(r)N_k''(r)\Big)\right.\\
&\qquad\qquad \left.+p(p-1)r^{p-2} \sum_{k\in vs(y)\setminus vs(x)} \sigma^p(y(k))\right)\,dr
\end{align*}
for $\alpha\ge\beta$. It follows that $\tau'$ is absolutely continuous and that (\ref{f.taudouble}) holds.
With the aid of Fubini's theorem, it also follows that (\ref{f.tauprime}) holds.

Since $\tau(\alpha)=0$ either for all $\alpha$ or for at most one $\alpha\in\RR$, it follows that $N'$
and $N''$ exists almost everywhere and that (\ref{f.tauprimeN}) and (\ref{f.taudoubleN}) hold for
almost every $\alpha\in\RR$.
\end{proof}

\subsection{The unmixed case}
If $X$ reflects disjointness and $p\in (2,\infty)$, then
the following assertion is true for $\ell_p(X)$.
\begin{proposition}\label{prop:4.1}
If $p\in (2,\infty)$ and $X$ is a $C^2$ smooth Banach space,
with a $1$-unconditional basis, such that $X$  reflects disjointness and
the first and second derivatives of the norm on $X$ are uniformly bounded on
$S_X$, then $\ell_p(X)$ reflects disjointness.
\end{proposition}
\begin{proof}
Let $\sigma$ denote the norm on $X$ and let $x,y\in \ell_p(X)$ with $x\not\in\linspan\{y\}$.
For each $k\in \mathbb{N}$ we define $N_k(\alpha)=\sigma(x(k)+\alpha y(k))$
and $\tau(\alpha)=N^p(\alpha)=\sum_{k=1}^\infty N_k^p(\alpha)$.
It follows from Theorem \ref{thm:2.2} that both $N$ and $\tau$ are continuously
differentiable on $\mathbb{R}$.
As $N$ is a convex function, the second derivative $N''(\alpha)$ exists almost
everywhere and the first condition in Definition \ref{def:3.1} is satisfied.

Due to Lemma \ref{lem.L1}, there exists a subset $A$ of $[0,1]$ with Lebesgue measure $1$
such that (\ref{f.tauprime})--(\ref{f.taudoubleN}) hold for all $\alpha\in A$.
Now assume that $N'(0)=0$ and that $N''(\alpha)$ converges to $0$, as
$\alpha \to 0$ along a subset of $[0,1]$ of full measure.
As $N'$ is continuous near $0$, it follows that $\tau'(0)=0$ and
$\tau''(\alpha)$ converges to $0$, as $\alpha \to 0$ along a subset of
$[0,1]$ of full measure.
Since each term in the sums in (\ref{f.taudouble}) is nonnegative we deduce for each
$k\in vs(x)\cap vs(y)$, that $N''_k(\alpha)$ converges to $0$, as $\alpha \to 0$
along a subset of $[0,1]$ of full measure and the continuity of $N_k'$ implies that
$N'_k(0)=0$. If $X$ reflects disjointness, we find that $x(k)$ and $y(k)$ are disjoint 
in $X$ for all $k \in vs(x)\cap vs(y)$.
Thus, $x$ and $y$ are disjoint in $\ell_p(X)$ and hence the second condition in
Definition \ref{def:3.1} is satisfied.

To prove the third condition, we assume that $x$ and $y$ are disjoint and $y$ is simple.
As $y$ is simple, the sums in (\ref{f.tauprime}) and (\ref{f.taudouble}) consist of finitely many terms.
Since $p>2$ and $X$ reflects disjointness, we find that $\tau'(0)=0$ and $\tau''(\alpha)$ converges
to 0, as $\alpha\to 0$ along a subset of $[0,1]$ of full measure. By subsequently
using (\ref{f.tauprimeN}) and (\ref{f.taudoubleN}), we see that $N'(0)=0$ and $N''(\alpha)$ also converges
to 0, as $\alpha\to 0$ along a subset of $[0,1]$ of full measure, and we are done.
\end{proof}
A combination of Proposition \ref{prop:4.1} with Theorems \ref{thm:2.2} and
\ref{thm:3.4} immediately gives the following corollary.
\begin{theorem}\label{thm:4.2}
If $p\in (2,\infty)$ and $X$ is a $C^2$ smooth Banach space with a
$1$-unconditional basis, such that $X$ reflects disjointness and
the first and second derivatives of the norm on $X$ are uniformly bounded on $S_X$,
then every contractive projection on $\ell_p(X)$ is an
averaging operator and its range admits a block basis.
\end{theorem}
It is well-known that for each $q>2$ the space $\ell_q$ is $C^2$ smooth
and the first and second derivatives of the norm are uniformly bounded on $S_X$
(see e.g. \cite[Chapter~V]{DGZ}).
Furthermore, since $\RR$ reflects disjointness, it follows from Proposition \ref{prop:4.1}
that $\ell_q$ reflects disjointness, if $q>2$.
Therefore we have the following result.
\begin{corollary}\label{cor:4.3}
If $p,q\in (2,\infty)$ or $p,q\in (1,2)$, then the range of every
contractive projection on $\ell_p(\ell_q)$ has a block basis.
\end{corollary}
\begin{proof}
The case $p$ and $q$ in $(2,\infty)$ is an immediate consequence of
Theorem \ref{thm:4.2}.
For $p,q\in (1,2)$ the assertion follows from the fact that the dual of
$\ell_p(\ell_q)$ is equal to $\ell_{p^*}(\ell_{q^*})$, where $1/p+1/p^*=1$ and
$1/q+1/q^*=1$.
\end{proof}

\subsection{The mixed case}
In the mixed case the space $\ell_p(\ell_q)$ is not $C^2$ smooth and it does not
reflect disjointness.
Therefore we can not apply Theorem \ref{thm:3.4} to show that every contractive
projection is an averaging operator.
Instead of using Theorem \ref{thm:3.4} we show that every contractive projection on
$\ell_p(\ell_q)$ is semi-band preserving.
The argument is quite involved and split up into several steps.
We begin by proving the following proposition.
\begin{proposition} \label{prop:4.4}
If $p\in (1,2)$ and $X$ is a $C^2$ smooth Banach space, with a $1$-unconditional basis,
such that $X$ reflects disjointness and the derivative of the norm $\sigma$ on $X$ is uniformly
bounded on $S_X$, then for each $x,y\in \ell_p(X)$ with $x\not\in\linspan\{y\}$ we have that
\begin{enumerate}[(a)]
\item the function $N(\alpha)=\|x+\alpha y\|$ is continuously differentiable and
$N''(\alpha)$ exists  almost everywhere;
\item if $vs(y)\subset vs(x)$, $N'(0)=0$, and $N''(\alpha)$
converges to $0$, as $\alpha\to 0$ along a subset of $[0,1]$ of full measure, then $x$ and $y$ are
disjoint;
\item if $x$ and $y$ are disjoint, $y$ is simple, and $vs(y)\subset vs(x)$, then $N'(0)=0$ and
$N''(\alpha)$ converges to $0$, as $\alpha\to 0$ along a subset of $[0,1]$ of full
measure;
\item if the second derivative of the norm of $X$ is uniformly bounded on $S_X$, if 
$vs(y)\subset vs(x)$, and $2\sigma(y(k))<\sigma(x(k))$ for all $k\in vs(x)$,
then there exists $A\subset [0,1]$ of full measure and $C>0$ such that
$N''(\alpha)\leq C$ for all $\alpha \in A$;
\item if $vs(y)\not\subset vs(x)$, then $N''(\alpha)\to \infty$, as $\alpha\to 0$
along a subset of $[0,1]$ of full measure.
\end{enumerate}
\end{proposition}
\begin{proof}
Let $\sigma$ denote the norm on $X$.
As in the proof of Proposition \ref{prop:4.1} we define
$N_k(\alpha)=\sigma(x(k)+\alpha y(k))$ for all $k\in \mathbb{N}$ and
$\tau(\alpha)= N^p(\alpha)=\sum_{k=1}^\infty N_k(\alpha)$.
By Theorem \ref{thm:2.2} both $\tau$ and $N$ are continuously differentiable on
$\mathbb{R}$.
Moreover, $N''$ exists almost everywhere, as $N$ is a convex function and hence
part (a) is satisfied.

Next, note that equations (\ref{f.tauprime})--(\ref{f.taudoubleN}) hold for almost every $\alpha\in\RR$.
Suppose that $vs(y)\subset vs(x)$, $N'(0)=0$, and
$N''(\alpha)$ converges to $0$, as $\alpha\to 0$, along a subset of $[0,1]$
of full measure.
Clearly the second sums in (\ref{f.tauprime}) and (\ref{f.taudouble}) are zero in that case.
As each term in the first sum of (\ref{f.tauprime}) and (\ref{f.taudouble}) is nonnegative,
we conclude from (\ref{f.tauprimeN}) and (\ref{f.taudoubleN}) that $N_k'(0)=0$, and
$N_k''(\alpha)$ converges to $0$, as $\alpha\to 0$, along a subset of $[0,1]$
of full measure for each $k\in vs(x)\cap vs(y)$.
Since $X$ reflects disjointness and $vs(y)\subset vs(x)$, it follows that
$x$ and $y$ are disjoint, which proves part (b).

To prove (c), note that if $x$ and $y$ are disjoint, $y$ is simple and $vs(y)\subset vs(x)$,
then for each $k\in vs(y)$ we have that $N_k'(0)=0$ and $N_k''(\alpha)$
converges to $0$, as $\alpha\to 0$, along a subset of $[0,1]$ of full measure.
As $y$ is simple, the sums in (\ref{f.tauprime}) and (\ref{f.taudouble}) consist
of finitely many terms, so that $\tau'(0)=0$ and $\tau''(\alpha)$ converges to $0$,
as $\alpha\to 0$, along a subset of $[0,1]$ of full measure.
The assertion now follows from equations (\ref{f.tauprimeN}) and (\ref{f.taudoubleN}).

To prove (d) assume that $vs(y)\subset vs(x)$.
Then the second sum in (\ref{f.taudouble}) vanishes.
For each $k\in vs(x)$  and $\alpha\in [0,1]$ we have that
\[
N_k(\alpha)\geq \sigma(x(k))-\alpha\sigma(y(k))
      > (1-\alpha/2)\sigma(x(k)) > \sigma(x(k))/2
\]
and
\[
N'_k(\alpha)= (x(k)+\alpha y(k))^*y(k)\leq \sigma(y(k))\leq \sigma(x(k))/2,
\]
where $(x(k)+\alpha y(k))^*$ denotes the norming functional of
$x(k)+\alpha y(k)$. This implies that $\tau'$ is bounded on $A$.

Since $D\sigma$ and $D^2\sigma$ are uniformly bounded on $S_X=\{x\in X: \sigma(x)=1\}$,
there exists a constant $c\in\mathbb{R}$ such that
$\| D\sigma(z)\|\leq c$ and $\|\sigma(z)D^2\sigma(z)\|\leq c$ 
for all $z\in X\setminus\{0\}$.
Thus, there exists $C>0$ such that
\begin{eqnarray*}
\lefteqn{|p(p-1)N_k^{p-2}(\alpha)N_k'(\alpha)^2 + pN_k^{p-1}(\alpha)N_k''(\alpha)| \leq } \\
   & & p(p-1)(\sigma(x(k)+\alpha y(k))^{p-2}c^2\sigma(y(k))^2) \\
   & &  \mbox{}  + p(\sigma(x(k)+\alpha y(k)))^{p-2}c\sigma(y(k))^2 \\
   & \leq & \frac{p(p-1)}{2}\sigma(x(k))^{p-2} c^2\sigma(x(k))^2 + \frac{p}{2}\sigma(x(k))^{p-2}
            c\sigma(x(k))^2\\
   & \leq & C\sigma(x(k))^p,
\end{eqnarray*}
as $p<2$. Therefore $\tau''$ is also bounded on $A$. It is now straightforward to deduce from
(\ref{f.tauprimeN}) and (\ref{f.taudoubleN}) that $N''$ is bounded on $A$.

Finally to prove (e) we assume that $vs(y)\not\subset vs(x)$.
In that case the second sum in (\ref{f.taudouble}) becomes unbounded, as
$\alpha \to 0$.
As every term in the first sum of (\ref{f.taudouble}) is nonnegative, we conclude
that $\tau''(\alpha)\to \infty$,
as $\alpha\to 0$, along a subset of $[0,1]$ on which $N''_k(\alpha)$ exists
for all $k\in vs(y)$. This subset of $[0,1]$ may be chosen such that it has full measure, as $vs(y)$ is countable.
By using (\ref{f.taudoubleN}) we deduce (e).
\end{proof}
To prove that a contractive projection on $\ell_p(X)$, with $p\in (1,2)$,
is semi-band preserving, we need to show  that $Px$ and $Py$ are disjoint,
whenever $Px$ and $y$, with $y$ simple, are disjoint.
To establish this, it is convenient to write
$y=y^1+y^2+y^3$, where $vs(y^1)\subset vs(Px)$, $vs(y^2)\cap vs(Pz)=\emptyset$
for all $z\in \ell_p(X)$, and $vs(y^3)\cap vs(Px)=\emptyset$, but $vs(y^3)\subset
vs(Pz)$ for some $z\in \ell_p(X)$.
The idea is to prove disjointness of $Py^i$ and $Px$ for $i=1,2$ and $3$
separately.
Let us start by analysing $y^1$.
\begin{lemma}\label{lem:4.5}
Let $p\in (1,2)$ and $X$ be a $C^2$ smooth Banach space, with a
$1$-unconditional basis, such that $X$ reflects disjointness and the
second derivative of the norm on $X$ is uniformly bounded on $S_X$.
If $P\colon \ell_p(X)\to \ell_p(X)$ is a contractive projection, then for each
$x,y\in \ell_p(X)$, with $vs(y)\subset vs(Px)$ and $y$ simple or $\sigma(y(k))\leq \sigma(Px(k))$ 
for all $k$, we have that $vs(Py)\subset vs(Px)$.
\end{lemma}
\begin{proof}
Let $x,y\in\ell_p(X)$ with $Px\neq 0$. For $\alpha\in\mathbb{R}$ define $\phi(\alpha)=\|Px+\alpha Py\|$ and
$\psi(\alpha)=\|Px+\alpha y\|$.
As $P$ is a contractive projection, $\phi(\alpha)\leq \psi(\alpha)$
for all $\alpha \in \mathbb{R}$ and $\phi(0)=\psi(0)$.
By rescaling $y$ we may assume without loss of generality that
$2\sigma(y(k))< \sigma(Px(k))$ for all $k\in vs(x)$.
Since $vs(y)\subset vs(Px)$,
we know by Proposition \ref{prop:4.4}(d) that there exists $C>0$
such that
\[
\{\alpha\in [0,1]\colon \psi''(\alpha) \mbox{ exists and }\psi''(\alpha)\leq C\}
\]
has full measure in $[0,1]$.
Hence it follows from Lemma \ref{lem:3.5}(c) that $\phi''(\alpha)$
does not go to infinity, as $\alpha\to 0$ along any subset of $[0,1]$
of full measure.
By using Proposition \ref{prop:4.4}(e) we conclude that $vs(Py)\subset vs(Px)$.
\end{proof}
This lemma has the following consequence.
\begin{lemma}\label{cor:4.6}
Let $p\in (1,2)$ and $X$ be a $C^2$ smooth Banach space, with a
$1$-unconditional basis, such that $X$ reflects disjointness and the second
derivative of the norm on $X$ is uniformly bounded on $S_X$.
If $P\colon \ell_p(X)\to \ell_p(X)$ is a contractive projection and
$x,y\in \ell_p(X)$, with $y$ simple, are such that $y$ and $Px$ are disjoint and
$vs(y)\subset vs(Px)$, then $Py$ and $Px$ are disjoint.
\end{lemma}
\begin{proof}
For $\alpha\in\RR$ define $\phi(\alpha)=\|Px+\alpha Py\|$ and $\psi(\alpha)=\|Px+\alpha y\|$.
As $P$ is a contractive projection, $\phi(\alpha)\le\psi(\alpha)$ for all $\alpha\in\RR$ and 
$\phi(0)=\psi(0)$.
We may assume that $Px\neq 0$. Due to Lemma \ref{lem.L1} and (\ref{f.tauprimeN}), $\phi$ and 
$\psi$ are $C^1$
and $\phi'$ and $\psi'$ are absolutely continuous functions on a neighborhood of $0$. As $y$ 
and $Px$ are disjoint, $y$ is simple, and $vs(y)\subset vs(Px)$, 
Proposition \ref{prop:4.4}(c) yields that $\psi'(0)=0$
and $\psi''(\alpha)\to 0$ as $\alpha\to 0$ along a full subset $A$ of $[0,1]$. By Lemma \ref{lem:3.5},
there exists a measurable set $E\subset[0,1]$ such that $E\cap (0,\delta)$ has positive measure for
all $\delta>0$ and such that $0\le \phi''(\alpha)\le\psi''(\alpha)$ for all $\alpha\in E$.
We may intersect $E$ with the full set $A$ and thus assume that $E\subset A$.
Then $\phi''(\alpha)\to 0$ as $\alpha\to 0$ along $E$. From (\ref{f.taudouble}) we infer that
$vs(Py)\subset vs(Px)$ and find for each $k\in vs(Px)\cap vs(Py)$ that,
\[
N_k^{p-1}(\alpha)N_k''(\alpha)\to 0\mbox{ as }\alpha\to 0\mbox{ along }E.
\]
Therefore, for $k\in vs(Py)$, $N_k''(\alpha)\to 0$ along $E$. Therefore $N_k''(\alpha)\to 0$ as $\alpha\to 0$
along a full subset of $[0,1]$, as $N_k''$ is continuous near $0$. Lemma \ref{lem:3.5} further
yields that $\phi'(0)=0$. From (\ref{f.tauprime}) it follows that $0\le N_k^{p-1}(\alpha)N_k'(\alpha)\le\phi'(\alpha)$
for almost every $\alpha$ and each $k\in vs(Py)$, so that the continuity of $\phi'$ and $N_k'$ gives $N_k'(0)=0$.
Since $X$ reflects disjointness, we obtain that $Px(k)$ and $Py(k)$ are disjoint for all $k\in vs(Py)$,  
and hence $Px$ and $Py$ are disjoint.
\end{proof}
A combination of Lemma \ref{lem:4.5} and \ref{cor:4.6} shows that $Py^1$ and $Px$ are disjoint. 
To prove disjointness for $y^2$, we shall use the following lemma.
\begin{lemma}\label{lem:4.7}
Let $p\in (1,2)$ and $X$ be a Banach space, with a
$1$-unconditional basis. 
If $P\colon \ell_p(X)\to \ell_p(X)$ is a contractive projection and
$y\in \ell_p(X)$, with $vs(y)\cap vs(Py) = \emptyset$ disjoint,
then $Py=0$.
\end{lemma}
\begin{proof}
Since $vs(y)\cap vs(Py)=\emptyset$, we have that
$\|Py +\alpha y\|^p=\|Py\|^p+\alpha^p\|y\|^p$
and $\|Py +\alpha Py\|^p= (1+\alpha)^p\|Py\|^p$ for all $\alpha\in [0,1]$.
As $P$ is a contractive projection, we deduce that
$(1+\alpha)^p\|Py\|^p\leq \|Py\|^p+\alpha^p\|y\|^p$,
so that
\[
\frac{\|Py\|^p}{\|y\|^p}\le\frac{\alpha^p}{(1+\alpha)^p-1}
\mbox{\quad for all }\alpha \in [0,1].
\]
Now note that as $p>1$,
\[
\lim_{\alpha\to 0} \frac{\alpha^p}{(1+\alpha)^p-1}=
\lim_{\alpha\to 0} \frac{p\alpha^{p-1}}{p(1+\alpha)^{p-1}}=0
\]
and hence $\|Py\|=0$.
\end{proof}
To prove disjointness for $y^3$ we need the following result.
\begin{lemma}\label{lem:4.8}
Let $p\in (1,2)$ and $X$ be a $C^2$ smooth Banach space, with a
$1$-unconditional basis, such that $X$ reflects disjointness and the
derivative of the norm on $X$ is uniformly bounded on $S_X$.
If $P\colon \ell_p(X)\to \ell_p(X)$ is a contractive projection
and $x,y,z\in\ell_p(X)$, with $y$ simple,  are such that
$vs(y)\subset vs(Pz)$ and $vs(y)\cap vs(Px)=\emptyset$, then there
exists $z'\in\ell_p(X)$ such that $vs(y)\subset vs(Pz')$ and
$vs(Pz')\cap vs(Px)=\emptyset$.
\end{lemma}
Before proving this lemma we give an auxiliary result.
\begin{lemma}\label{lem:4.9}
Let $p\in (1,2)$ and $X$ be a $C^2$ smooth Banach space, with a
$1$-unconditional basis, such that $X$ reflects disjointness and the
derivative of the norm on $X$ is uniformly bounded on $S_X$.
Suppose that $P\colon \ell_p(X)\to \ell_p(X)$ is a contractive projection
and denote
\[
\Sigma_P=\{A\subset\mathbb{N}\colon vs(Pu)=A\mbox{ for some }u\in \ell_p(X)\}.
\]
Then the following assertions are true:
\begin{enumerate}[(a)]
\item If $(A_i)_{i\in \mathbb{N}}\subset \Sigma_P$ and $A_1\supset A_2\supset\ldots\,$,
then $\cap_{i\in \mathbb{N}} A_i\in \Sigma_P$.
\item If $A,B\in \Sigma_P$ and $a\in A\setminus B$, then there exists
$D_a\in \Sigma_P$ such that $a\in D_a$ and $D_a\subset A\setminus B$.
\end{enumerate}
\end{lemma}
\begin{proof}
Let $(A_i)_{i\in \mathbb{N}}\subset \Sigma_P$
be such that $A_i\supset A_{i+1}$ for all $i\in \mathbb{N}$.
By definition there exist $u^i\in \ell_p(X)$ such that $vs(Pu^i)= A_i$
for each $i\in\mathbb{N}$.
Put $A=\cap_{i\in \mathbb{N}} A_i$ and let $w=(Pu^1)\chi_A$, where
$\chi_A$ is the indicator function of $A$.
Clearly, $vs(w)=A\subset vs(Pu^i)$ for all $i\in\mathbb{N}$.
Thus, Lemma \ref{lem:4.5} implies that $vs(Pw)\subset vs(Pu^i)$ for all
$i\in \mathbb{N}$ and hence $vs(Pw)\subset A$.
Put $B=vs(Pw)$ and remark that $(Pw)\chi_{A_1\setminus B}=0$.
Moreover,
\[
Pu^1= (Pu^1)\chi_{A_1\setminus A}+ (Pu^1)\chi_{A}
    = (Pu^1)\chi_{A_1\setminus A}+ w,
\]
so that $Pu^1=P^2u^1= P((Pu^1)\chi_{A_1\setminus A}) +Pw$.
Therefore
\[
(Pu^1)\chi_{A_1\setminus B}= P((Pu^1)\chi_{A_1\setminus A})\chi_{A_1\setminus B}+
                             (Pw)\chi_{A_1\setminus B}
    = P((Pu^1)\chi_{A_1\setminus A})\chi_{A_1\setminus B}.
\]
As $P$ is contractive, we find that
\begin{eqnarray*}
\|(Pu^1)\chi_{A_1\setminus B}\|
            &  =  & \|P((Pu^1)\chi_{A_1\setminus A})\chi_{A_1\setminus B}\| \\
            &\leq & \|P((Pu^1)\chi_{A_1\setminus A})\| \\
            &\leq & \|(Pu^1)\chi_{A_1\setminus A}\|.
\end{eqnarray*}
Since $vs(Pu^1)=A_1$ and $B\subset A$, we conclude that $A=B=vs(Pw)$ and hence
$A\in \Sigma_P$.

To prove the second assertion let $A, B\in \Sigma_P$ and $a\in A \setminus B$.
If $A\cap B=\emptyset$, then we can take $D_a=A$.
So, suppose that $A\cap B=B_0$ is not empty.
By the first assertion, $B_0\in\Sigma$.
Now let $u\in\ell_p(X)$ be in the range of $P$ and $vs(u)=A$.
It follows from Lemma \ref{lem:4.5} that
\[
vs(P(u\chi_{A\setminus B_0}))\subset A\mbox{\quad and \quad }
vs(P(u\chi_{B_0}))\subset B_0.
\]
But also, $P(u\chi_{A\setminus B_0}) + P(u\chi_{B_0}) = Pu = u$, 
so that
\begin{equation}\label{eq:4.8}
P(u\chi_{A\setminus B_0})\chi_{A\setminus B_0}= u\chi_{A\setminus B_0}.
\end{equation}
As
\[
P(u\chi_{A\setminus B_0}) =
P(u\chi_{A\setminus B_0})\chi_{A\setminus B_0}+ P(u\chi_{A\setminus B_0})\chi_{B_0}
\]
and $\|P(u\chi_{A\setminus B_0})\|\leq \| u\chi_{A\setminus B_0}\|$, it follows
from (\ref{eq:4.8}) that $\|P(u\chi_{A\setminus  B_0})\chi_{B_0}\|=0$.
Thus, $P(u\chi_{A\setminus B_0})\chi_{A\setminus B_0 }= Pu\chi_{A\setminus B_0}$
and we can take $D_a = A\setminus B_0$.
\end{proof}
Using this lemma it is now straightforward to prove Lemma \ref{lem:4.8}.
\begin{proof}[Proof of Lemma \ref{lem:4.8}]
Let $y$ be  simple  and for each $k\in vs(y)$, let $D_k$ be a set in $\Sigma_P$ given in
Lemma \ref{lem:4.9}(b), where $A=vs(Pz)$ and $B=vs(Px)$.
Put $D=\cup_{k\in vs(y)} D_k$ and note that, as $vs(y)$ is finite, there exists
$z'\in \ell_p(X)$ such that $vs(Pz')=D$ and this completes the proof.
\end{proof}
A combination of the lemmas now yields the following theorem.
\begin{theorem}\label{thm:4.10}
If $p\in (1,2)$ and $X$ is a $C^2$ smooth Banach space with a
$1$-unconditional basis, such that $X$ reflects disjointness and
the second derivative of the norm on $X$ is uniformly bounded on $S_X$,
then every contractive projection on $\ell_p(X)$ is an
averaging operator and its range admits a block basis.
\end{theorem}
\begin{proof}
Let $x,y\in \ell_p(X)$, with $y$ simple, be such that $Px$ and $y$ are disjoint.
Write $y=y^1+y^2+y^3$ with each $y^i$ simple, where $vs(y^1)\subset vs(Px)$, $vs(y^2)\cap vs(Pw)=\emptyset$
for all $w\in \ell_p(X)$, and $vs(y^3)\subset vs(Pz)$ for some $z\in \ell_p(X)$, with
$vs(Pz)\cap vs(Px)=\emptyset$.
Then it follows from Corollary \ref{cor:4.6} that $Py^1$ and $Px$ are disjoint.
Moreover, Lemma \ref{lem:4.7} implies that $Py^2=0$ and hence $Py^2$ and $Px$ are
disjoint.
For $y^3$ we find by Lemma \ref{lem:4.8} that there exists $z'\in\ell_p(X)$ such that
$vs(y^3)\subset vs(Pz')$ and $vs(Pz')\cap vs(Px)=\emptyset$.
In addition, it follows from Lemma \ref{lem:4.5} that $vs(Py^3)\subset vs(Pz')$, 
so that $Py^3$ and $Px$ are disjoint.
\end{proof}
Theorem \ref{thm:4.10} has the following consequence for $\ell_p(\ell_q)$ spaces. 
\begin{corollary}\label{cor:4.11}
If $p\in (1,2)$ and $q\in (2,\infty)$, or, the other way around,
then the range of every contractive projection on $\ell_p(\ell_q)$
has a block basis.
\end{corollary}
\begin{proof}
The proof follows from Theorem \ref{thm:4.10}, the fact that the dual of
$\ell_p(\ell_q)$ is equal to $\ell_{p^*}(\ell_{q^*})$, where $1/p+1/p^*=1$ and
$1/q+1/q^*=1$, and the fact that $\ell_r$ reflects disjointness if $r>2$.
\end{proof}

\section{Conclusions}
Combining Corollaries \ref{cor:4.3} and \ref{cor:4.11} with the results from \cite[Section 5]{Ran1} 
yields the following characterization  of 1-complemented subspaces of $\ell_p(\ell_q)$ spaces. 
\begin{theorem}\label{thm:5.1}
If $p,q\in (1,\infty)$, with $p,q\neq 2$, and $Y$ is a subspace of $\ell_p(\ell_q)$,
then $Y$ is the range of a contractive projection on $\ell_p(\ell_q)$ if and only
if there exists a basis $\{y^i\}_{i=1}^{\dim Y}$ for $Y$ such that  for each $i\neq j$
either $vs(y^i)\cap vs(y^j)=\emptyset$ or $vs(y^i)=vs(y^j)$ and in that case
and $y^i(k)$ and $y^j(k)$ are disjoint and $\|y^i(k)\|_q=\|y^j(k)\|_q$ for all $k\in vs(y^i)$.
\end{theorem}
Of course, it would be interesting to see if this theorem can be extended to general 
vector-valued $L_p(L_q)$ spaces. 

Concerning Theorem \ref{thm:3.4}, we like to remark that if $X$ is a Banach sequence space 
with norm, $\|x\|=\|x\|_p+\|x\|_q$, and $p,q>2$, then $X$ reflects disjointness and hence every
contractive projection on $X$ is an averaging operator and its range has a block basis.
However, the theorem can not be applied if $p=2$ and $q>2$.
Nevertheless we believe that the same assertion is true, but, as yet, 
we can not prove it.
In connection with this problem it is worth mentioning a general conjecture of Randrianantoanina
\cite[Conjecture 7.9]{Ran5}, which asserts that if $X$ is strictly monotone
Banach sequence space, with a $1$-unconditional basis, and $X$ does not contain an isometric copy of a 
Euclidean plane, then the range of every contractive projection on $X$ admits a block basis.
It is known \cite{Ran1} that this conjecture is true in complex Banach spaces.

\begin{small}
\vspace{4mm} 
\noindent
\begin{tabular}{lll}
Bas Lemmens          & \mbox{}  & Beata Randrianantoanina      \\
Mathematics Institute & \mbox{}  & Department of Mathematics and Statistics\\
University of Warwick & \mbox{} & Miami University    \\
CV4 7AL Coventry & \mbox{} & Oxford, OH 45056 \\
United Kingdom & \mbox{} & U.S.A. \\
lemmens@maths.warwick.ac.uk & \mbox{} & randrib@muohio.edu\\
\end{tabular}

\vspace{6mm}\noindent
\begin{tabular}{l}
Onno van Gaans\\
Mathematical Insitute\\
Leiden University\\
P.O. Box 9512\\
2300 RA Leiden \\
The Netherlands\\
vangaans@math.leidenuniv.nl\\
\end{tabular}
\end{small}

\begin{thebibliography}{10}
\bibitem{An} T. Ando, Contractive projections in ${L}_p$ spaces.
\emph{Pacific J. Math.} \textbf{17} (1966), 391--405.

\bibitem{BL} S.J. Bernau and H.E. Lacey,
The range of a contractive projection on an ${L}_p$-space.
\emph{Pacific J. Math.} \textbf{53} (1974), 21--41.

\bibitem{Bo} F. Bohnenblust, Subspaces of $\ell_{p,n}$ spaces.
\emph{Amer. J. Math.} \textbf{63} (1941), 64--72.

\bibitem{BHL} B. Bru, H. Heinich, J-C. Lootgieter,
It\'eration des applications de dualit\'e dans les espaces d'Orlicz.
\emph{C. R. Acad. Sci. Paris S\'er. I Math.} \textbf{303}(15) (1986), 745--747.

\bibitem{CP} K.W. Cheney and K.H. Price,
Minimal projections. \emph{1970 Approximation Theory (Proc. Sympos.,
Lancaster, 1969)} 261--289, Academic Press, London.

\bibitem{DGZ}
R.~Deville, G.~Godefroy, and V.~Zizler,
\newblock {\em Smoothness and renormings in {B}anach spaces}, volume~64 of {\em
  Pitman Monographs and Surveys in Pure and Applied Mathematics}.
\newblock Longman Scientific \& Technical, Harlow, 1993.

\bibitem{DU}
J.~Diestel and J.~J. Uhl, Jr.
\newblock {\em Vector measures}.
\newblock American Mathematical Society, Providence, R.I., 1977.
\newblock With a foreword by B. J. Pettis, Mathematical Surveys, No. 15.

\bibitem{Dieu} J. Dieudonn\'e,
\emph{Foundations of modern analysis}.
Pure and Applied Mathematics, Vol. X Academic Press, New York-London 1960.

\bibitem{Do} R.G. Douglas,
Contractive projections on an ${L}_1$-space.
\emph{Pacific J. Math.} \textbf{15} (1965), 443--462.

\bibitem{Doust} I. Doust,
Contractive projections on Lebesgue-Bochner spaces.
\emph{Function spaces (Edwardsville, IL, 1994)},
101--109, Lecture Notes in Pure and Appl. Math., \textbf{172}, Dekker,
New York, 1995.

\bibitem{JKL} J. Jamison and A. Kami\'nska and G. Lewicki,
One-complemented subspaces of {M}usielak-{O}rlicz spaces.
\emph{J. Approx. Theory} \textbf{130}(1) (2004), 1--37.

\bibitem{KW} N.J. Kalton and G.V. Wood, Orthonormal systems in 
{B}anach spaces and their applications. 
\emph{Math. Proc. Cambridge Philos. Soc.} \textbf{79}(3) (1976), 493--510.


\bibitem{Ka} S. Kakutani, Some characterizations of {E}uclidean spaces.
\emph{Jap. J. Math.} \textbf{16} (1940), 93--97.

\bibitem{Kol} A. L. Koldobsky, Isometries of $L\sb p(X;L\sb q)$ and
equimeasurability.  \emph{Indiana Univ. Math. J.}
\textbf{40}(2) (1991), 677--705.

\bibitem{LvG} B. Lemmens and O. W. van Gaans,
On one-complemented subspaces of Minkowski spaces with smooth Riesz norms.
\emph{Rocky Mountain J. Math.}, to appear.

\bibitem{LS} I.E. Leonard and  K. Sundaresan,
Geometry of Lebesgue-Bochner function spaces -- smoothness.
\emph{Trans. Amer. Math. Soc.} \textbf{198} (1974), 229--251.

\bibitem{Le79} D.R. Lewis, Ellipsoids defined by {B}anach ideal norms.
\emph{Mathematika} \textbf{26}(1) (1979), 18--29. 

\bibitem{LT} J. Lindenstrauss and L. Tzafriri,
\emph{Classical Banach Spaces I. Sequence Spaces}, Springer-Verlag,
Berlin-New York, 1977.

\bibitem{Ran1} B. Randrianantoanina,
$1$-complemented subspaces of spaces with $1$-unconditional bases.
\emph{Canad. J. Math.} \textbf{49}(6) (1997), 1242--1264.

\bibitem {Ran5} B. Randrianantoanina,
Norm-one projections in Banach spaces.
\emph{Taiwanese J. Math.}  \textbf{5}(1) (2001), 35--95.

\bibitem{Ran6} B. Randrianantoanina,
Contractive projections in Orlicz sequence spaces.
\emph{Abstr. Appl. Anal.} \textbf{2004}(2) (2004), 133--145.

\bibitem{Ran7} B. Randrianantoanina,
A disjointness type property of conditional expectation operators.
\emph{Colloq. Math.} \textbf{102} (2005),  9--20.

\bibitem{Ray} Y. Raynaud,
The range of a contractive projection in $L\sb p(H)$.
\emph{Rev. Mat. Complut.} \textbf{17}(2) (2004), 485--512.

\bibitem{Su} K. Sundaresan, Smooth Banach spaces.
\emph{Math. Ann.} \textbf{173} (1967), 191--199.

\bibitem{SS} K. Sundaresan and S. Swaminathan,
\emph{Geometry and nonlinear analysis in Banach spaces}.
Lecture Notes in Mathematics, \textbf{1131}. Springer-Verlag, Berlin, 1985.

\bibitem{Tz} L. Tzafriri,
Remarks on contractive projections in $L\sb{p}$-spaces.
\emph{Israel J. Math.} \textbf{7} (1969),  9--15.

\end{thebibliography}
\end{document}